\documentclass{siamltex}
\usepackage{amsmath,amssymb,amscd,fancybox,ifthen,float,epsfig}
\usepackage{microtype}

\usepackage{tikz}
\usetikzlibrary{positioning}
\usetikzlibrary{calc}
\usetikzlibrary{chains}
\usetikzlibrary{shapes}

\usepackage{ifthen}

\numberwithin{equation}{section}
\floatstyle{plain}
\restylefloat{figure}

\newcommand\cK{\mathcal K}

\newcommand\tm{\widetilde m}

\newcommand\hvarphi{\widehat \varphi}

\newcommand\bn{\boldsymbol n}
\newcommand\btau{\boldsymbol \tau}

\newcommand\bz{\bar{z}}

\newcommand\talpha{\widetilde{\alpha}}
\newcommand\halpha{\widehat{\alpha}}

\newcommand\cC{\mathcal{C}}

\renewcommand\Re{\operatorname{Re}}
\renewcommand\Im{\operatorname{Im}}

\newcommand\bbC{\mathbb C}

\newcommand\bbR{\mathbb R}

\newcommand\pa{\partial}

\newcommand\restrictedto{\upharpoonright}

\newcommand\CI{{\mathcal C}^{\infty}}

\newtheorem{remark}{Remark}

\newcommand{\doi}[1]{\href{http://dx.doi.org/#1}{doi: #1}}

\def\dlmf#1{\href{http://dlmf.nist.gov/#1}{#1}}

\newboolean{includefigures}
\setboolean{includefigures}{true}

\usepackage{hyperref}

\begin{document}


\title{On the convergence of local expansions of layer~potentials}

\author{Charles L. Epstein\thanks{
    Depts. of Mathematics and Radiology, University of Pennsylvania,
    209 South 33rd Street, Philadelphia, PA 19104. E-mail:
    {cle@math.upenn.edu}.
    Research partially supported by NSF grants
    DMS09-35165 and DMS12-05851.} \and
 Leslie Greengard\thanks{Courant Institute,
    New York University, 251 Mercer Street, New York, NY 10012.
    E-mails: {\{greengard,kloeckner\}@cims.nyu.edu}.
    Research partially supported by the U.S. Department of Energy under
    contract DEFG0288ER25053 and by the Air Force Office of Scientific Research
    under NSSEFF Program Award FA9550-10-1-0180.} \and
    Andreas Kl\"ockner\footnotemark[2]
}
\date{}

\maketitle

\begin{abstract}
In a recently developed quadrature method (quadrature by expansion or QBX),
it was demonstrated that weakly singular or singular layer potentials can be evaluated
rapidly and accurately {\em on surface}
by making use of local expansions about carefully chosen off-surface points.
In this paper, we derive estimates for the rate of convergence of these local
expansions, providing the analytic foundation for the QBX method. The estimates
may also be of mathematical interest, particularly for microlocal or
asymptotic analysis in potential theory.
\end{abstract}

\begin{keywords}
Integral equations, layer potential, quadrature, singular integrals,
spherical harmonics, expansion, Laplace equation, Helmholtz equation.
\end{keywords}

\begin{AMS}
  65R20, 65N38, 65N80, 31A10, 65D32
\end{AMS}

\pagestyle{myheadings}
\thispagestyle{plain}
\markboth{C. EPSTEIN, L. GREENGARD AND A. KL\"OCKNER}{ON THE CONVERGENCE OF LOCAL EXPANSIONS OF LAYER~POTENTIALS}


\section{Introduction}
The paper~\cite{QBX} describes a new method for the evaluation of layer
potentials referred to as `quadrature by expansion'. This method, denoted
by QBX, is a technique for evaluating layer potentials \emph{on surface},
which has simplified the development of fast, high-order accurate solvers for
boundary integral equations. It is particularly easy to implement because
it does not require the evaluation of integrals with non-smooth integrands.
Here, we present the analytic foundations for QBX, which consist mainly in the
establishment of decay estimates for one-sided local expansions induced by the
layer potentials themselves.  While introduced in the present context for the purpose
of numerical analysis, these estimates may be of interest in
their own right in asymptotic analysis and potential theory.

For simplicity we assume that $\Gamma\subset\bbR^n$ is a
smooth compact hypersurface separating $\bbR^n$ into two components
$\bbR^n\setminus\Gamma= D_+\cup D_-.$
Throughout the paper, we let $D_-$ denote the bounded
component and $D_+$ the unbounded component.
By a layer potential, we mean an integral of the form
\begin{equation}
  F_{\pm}(x)=\int\limits_{\Gamma}k(x,y)\varphi(y)dS(y),\text{ for }x\in D_{\pm},
\label{genlp}
\end{equation}
where $\varphi(y)$ is a smooth density on $\Gamma$ and
$k(x,y)$ is the Schwartz kernel, defined on
$\bbR^n\times\bbR^n,$ for a pseudodifferential operator, $\cK,$ which
satisfies the transmission condition, see Section 18.2 in~\cite{hormander3}. We are primarily interested in
the case where $k(x,y)$ is the Green's function
for the Laplace or Helmholtz equation (or the Cauchy kernel), but more
general cases also arise in a number of applications.

The kernels $k(x,y)$ in these calculations have
singularities on the diagonal (when $x=y$) which have hampered the development of
stable, efficient, and high order accurate methods for their evaluation,
particularly on surfaces in $\bbR^3$.
The most frequently used high-order methods to date that are suitable for
use with patch-based surface discretizations have largely been
based on product integration. In this approach, the integral
of the kernel multiplied by a
piecewise polynomial approximation of the density $\varphi$
on a piecewise smooth approximation of the boundary $\Gamma$ is computed
to high order using a mixture of
analysis, linear algebraic techniques, and optimization.
These have been developed in both two and three dimensions for a variety
of kernels
\cite{bremer_nonlinear_2010,helsing_2008b,strain_1995,yarvin_generalized_1998}.
Other methods, more closely related to QBX, are based on regularization of the kernel,
combined with smooth quadrature rules and asymptotic correction
\cite{beale_lai_2001, goodman_1990, haroldson_1998, lowengrub_1993, schwab_1992}.
Remarkably little use has been made, however, of the fact that the functions
$F_{\pm}(x)$ have smooth extensions to $\overline{D_{+}}$ and
$\overline{D_{-}}$. The limits
from the two sides, of course, often differ.

\begin{remark}
An exception is \cite{delves_1967,ioakimidis,lyness_1967}, in which
the authors do make use of a local expansion
about an off-surface point. In those papers, however, the viewpoint is
global. In essence, a single expansion center is introduced
(for Cauchy integrals with analytic data), with a radius of convergence determined by
the location of the nearest singularity of the analytic function itself.

In our case, we consider boundary data of finite differentiability, which is
not the restriction of a real- or complex-analytic function and hence make no
assumptions about the location of the nearest underlying singularity.  We are
interested in the behavior of {\em local} expansions in balls centered at
points close to the boundary $\Gamma,$ with radius  equal to the distance from
the center to the boundary. Due to the nature of the kernel functions,
$k(x,y),$ we can establish error estimates, valid in the \emph{closed} ball, that
only depend on the smoothness of the data and boundary near to the
intersection of $\Gamma$ with the boundary of the ball.
\end{remark}

A more detailed review of existing approaches to quadrature is given in \cite{QBX}
and the textbooks \cite{atkinson_1997,brebbia_1984,kress_1999}.

\ifthenelse{\boolean{includefigures}}{
\begin{figure}[h]
 \begin{center}
    \begin{tikzpicture}
      \coordinate (c) at (0,0) ;
      \path (c) ++(45+75:2.75) coordinate (s);
      \path (c) ++(45+15:2) coordinate (t);
      \path (c) ++(45+15:2) coordinate (t-to-curve);

      \path (s) ++(-1.5,0) coordinate (curve-before);
      \path (t) ++(2,-0.7) coordinate (curve-after);
      \draw [thick]
        (curve-before)
        ..controls +(30:0.7) and +(180-45:0.7) ..
        (s)
        node [pos=0.2,anchor=north] {$\Gamma$}
        ..controls +(-45:0.7) and +(45+90+15:1.25) ..
        (t-to-curve)
        ..controls +(45+270+15:0.7) and +(190:0.7) ..
        (curve-after) ;

      \fill (c) circle (1pt);
      \fill (t) circle (1pt);
      \fill (s) circle (1pt);

      \draw [dotted] ++(-20:2) arc (-20:45+105:2);

      \node at (105:2) {$r$};
      \draw (c) -- (t) ;
      \draw (c) -- (s) ;
      \node at (c) [anchor=north west] {$x_c$};
      \node at (s) [anchor=south west] {$y$};
      \node at (t) [anchor=south west,fill=white,xshift=0.5mm,inner sep=1mm] {$x_0$};
      \draw [dashed,->] (c) ++(-1, 0) -- ++(3.5,0);

      \draw (c) ++(0.5,0) arc (0:45+75:0.5);
      \path (c) ++(45*0.8+75*0.8:0.5)
        node [anchor=south] {$\phi$};

      \draw (c) ++(0.75,0) arc (0:45+15:0.75) ;
      \path (c) ++(45/2+15/2:0.75)
        node [anchor=west] {$\theta$};

      \draw [very thin] (t) -- (t-to-curve);
      \draw [very thin]
        let
          \p1 = ($ 0.1*(t-to-curve) - 0.1*(c) $),
          \p2 = (-\y1,\x1)
        in
        ($(t-to-curve)!0.1!(c)$) -- ++(\p2) -- ++(\p1) ;
    \end{tikzpicture}
 \end{center}
  \caption{A QBX expansion near a curve $\Gamma$, with a source point
    $y$, a target point $x_0$, and an expansion center $x_c$.
    The angles $\phi$ and $\theta$ are mainly used in Section~\ref{sec:two-d-helmholtz}.}
\end{figure}
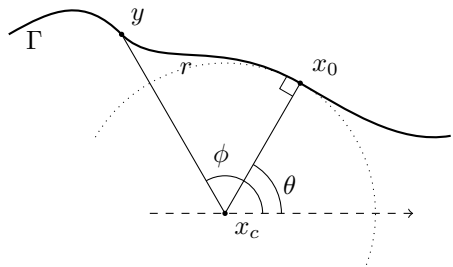
}{}

The QBX approach is based on the fact
that, because of its one-sided smoothness,
$F_{\pm}(x)$ supports a
Taylor series expansion about $x$. More precisely, suppose that
$x_0\in \Gamma$, that $x_c \in D_{-}$ and that there is a ball of radius
$r$ about $x_c$, $B_r(x_c)\subset D_{-}$,
such that $x_0\in bB_r(x_c)$, where $bB_r(x_c)$ denotes the boundary of the ball.
In other words, $x_0$ is a point of tangency of the ball and the surface $\Gamma$.
(The analysis for the $+$-case is essentially
identical.) The fact that $F_{-}$ is smooth in $\overline{B_r(x_c)}$
and the remainder theorem for Taylor series shows that
\begin{equation} \label{genexp}
  \left|F_{-}(x_0)-\sum_{\{\alpha:|\alpha|\leq k\}}
\frac{\pa_{x}^{\alpha}F_{-}(x_c)}{\alpha!}(x_c-x_0)^{\alpha}\right|\leq M_{k}(F_{-})r^{k+1}.
\end{equation}
Here $\alpha$ denotes a standard multi-index in $n$ dimensions, and
$M_{k}(F_{-})$ is a constant depending on the $\cC^{k+1}$-norm of $F_{-}$ in
$B_r(x_c).$ Because $\cK$ satisfies the transmission condition, this quantity
is bounded by the H\"older $\cC^{k+2+m,\alpha}$-norm of $\varphi,$ where the
$m$ is the order of $\cK,$ and $\alpha$ is any positive number. The two key
observations underlying QBX are that the error in the expansion \eqref{genexp}
can be controlled and that the derivatives $\{\pa_{x}^{\alpha} F_{-}(x_c)\}$
can be computed accurately using only standard quadrature techniques for {\em
  smooth} functions.

Briefly stated, QBX consists of a three-step procedure: For each boundary point $x_0$, we

\begin{enumerate}
\item find a nearby off-surface point $x_c$,
\item construct a local expansion of the form \eqref{genexp} about the center $x_c$
of the field induced by the layer potential, and
\item evaluate a partial sum of the local expansion from the previous step at $x_0$.
\end{enumerate}
For kernels that satisfy a classical PDE, more efficient series expansion
are preferable, based on standard separation of variables. For example, if
$\cK$ is a single or double layer potential for the Helmholtz
equation in three dimensions, for which the kernel is
$k(x,y)=e^{ik|x-y|}/{(4\pi|x-y|)},$
then it is convenient to use a (doubly-indexed) spherical harmonic expansion instead
of a (triply-indexed) Taylor series.

We  restrict our attention in this paper to kernels $k(x,y)$ which
are Green's functions of some constant-coefficient elliptic partial
differential operator, which are well-known to
have their singular support along the diagonal. The
analysis below, therefore, applies to rectifiable hypersurfaces $\Gamma,$
provided that $\Gamma$ is smooth ($\CI$) in a neighborhood of the point, $x,$
where we are trying to estimate the value of $F_{\pm}(x).$ To simplify the
discussion below we assume that $\Gamma$ is globally infinitely
differentiable. This is not necessary; from our analysis it is clear that the
accuracy of this approach near a point $x_0$ on the boundary depends only on the
smoothness of the boundary and the data near to $x_0$.

The purpose of this paper is to provide rigorous error estimates for the QBX
procedure outlined above. In Sections~\ref{sec:two-d} and \ref{sec:three-d}, we
derive decay estimates for analogues of $M_k(F_-)$ in \eqref{genexp} for a
selection of commonly encountered kernels. In Section~\ref{sec:quad}, we couple
this analysis with quadrature error estimates for the numerical computation of
$\{\pa_{x}^{\alpha} F_{-}(x_c)\}$. The corresponding analysis for the case of a
general kernel is omitted. The steps involved are essentially the same as for
the cases studied here.
\section{The two-dimensional case}
\label{sec:two-d}

We first consider two dimensional examples arising in potential and
scattering theory, namely the Cauchy kernel, the Green's function
for the Laplace equation, and the Green's function
for the Helmholtz equation.

\subsection{The Cauchy Kernel}
Let $D_-$ be a bounded
domain in $\bbC=\bbR^2$ with a smooth boundary $\Gamma=bD_-,$ and $\varphi$ a
function on $\Gamma.$ The Cauchy formula defines an analytic function
$f$ in $\bbC\setminus \Gamma:$
\begin{equation}\label{eqn3}
  f(z)=\frac{1}{2\pi i}\int\limits_{\Gamma}\frac{\varphi(\zeta)d\zeta}{\zeta-z}.
\end{equation}
We denote the restrictions $f_-=f\restrictedto_{D_-}$ and
$f_+=f\restrictedto_{D_+}$ and recall the well-known jump formula:
\begin{equation}
  [f_+-f_-]\restrictedto_{\Gamma}=\varphi.
\end{equation}

The crucial fact in our analysis is that the functions $f_{\pm}$ extend to the
closures of their domains of definition to be essentially as smooth as
$\varphi.$

\begin{remark}
  If smoothness is measured in terms of $L^2$-Sobolev norms,
  then $f_{\pm}$ have a half-derivative more than $\varphi,$ whereas,
  if we use standard H\"older norms, then $f_{\pm}$ have the same
  regularity, up to the boundary as the boundary data. Of course
  $f_{\pm}$ are analytic functions in the interior, so here we are
  primarily speaking of the boundary regularity.
\end{remark}

For purposes of definiteness, we suppose that $z_c\in D_-$, that
the disk $B_r(z_c)\subset D_-$, and that $z_0\in bB_r(z_c)\cap \Gamma$,
where $z_0 = z_c + re^{i\theta_0}$. In
this disk:
\begin{equation} \label{taylorz}
  f_-(z)=\sum_{j=0}^{\infty}\frac{f_-^{(j)}(z_c)}{j!}(z-z_c)^j,
\end{equation}
where $f^{(j)}(z)=\pa_z^jf(z).$
In general this power series converges in $B_r(z_c),$ but in no larger
disk. If $\varphi$ has more than a half an $L^2$-derivative, or is
H\"older continuous, then the series representation for $f_-$
converges uniformly on $\overline{B_r(z_c)}.$

In the QBX approach, we
would like to approximate the integral
$f_-$ from~\eqref{eqn3}, at the boundary point $z_0,$ by using a finite partial
sum of the Taylor series~\eqref{taylorz}, rather than a quadrature rule that
is designed to handle the singular kernel head-on.
Because of the special
structure of the Cauchy kernel, two different approaches are available
to estimate the error:
\begin{equation}
  e_N(z)= f_-(z)-\sum_{j=0}^N\frac{ f^{(j)}_-(z_c)}{j!}(z-z_c)^j,
\end{equation}
for $z=z_c+re^{i\theta_0}.$ The first method does not generalize, but
gives a simple explanation as to why this approach works, whereas the
second gives a clear path to an error estimate in the general case.

Let $w(t)$ be an arc-length parametrization for $\Gamma.$ The Cauchy
integral is then given by
\begin{equation}
  f(z)=\frac{1}{2\pi i}\int\limits_{0}^L\frac{\varphi(w(t))w'(t)}{w(t)-z}dt.
\end{equation}
To simplify notation, let us assume (without loss of generality) that $z_c = 0$.
For $z\in B_r(0)$, by considering the series for $(1-z/w(t))^{-1}$, we
obtain
\begin{equation}
  e_N(z)=\frac{1}{2\pi i}\sum_{j=N+1}^{\infty}\int\limits_{0}^L
\varphi(w(t))\left(\frac{z}{w(t)}\right)^j\frac{w'(t)}{w(t)}dt.
\end{equation}
As
\begin{equation}
  \frac{w'(t)}{(w(t))^{j+1}}=\frac{-1}{j}\pa_t\frac{1}{(w(t))^j},
\end{equation}
then we can integrate by parts to obtain
\begin{equation}
  e_N(z)=\frac{z}{2\pi i}\sum_{j=N+1}^{\infty}\int\limits_{0}^L
\pa_t[\varphi(w(t))]\left(\frac{z}{w(t)}\right)^{j-1}\frac{dt}{j\cdot w(t)}.
\end{equation}
In order to integrate by parts again, we use the fact that, for $j>1,$
\begin{equation}\label{eqn11.01}
  \frac{1}{[w(t)]^j}=\frac{-1}{(j-1)w'(t)}\pa_t\left( \frac{1}{[w(t)]^{j-1}}\right).
\end{equation}
If we define the first order differential operator $D_t$ by
\begin{equation}\label{eqn12.01}
  D_tg=\pa_t\left(\frac{g(t)}{w'(t)}\right),
\end{equation}
then we can  integrate by parts $N$ more times to obtain that
\begin{equation}
  e_N(z)=\frac{z^{N+1}}{2\pi i}
  \sum_{j=N+1}^{\infty}\int\limits_{0}^L
  D_t^N\pa_t[\varphi(w(t))]
  \bigg(\underbrace{\frac{z}{w(t)}}_{|\cdot|\le 1}\bigg)^{j-N-1}
  \frac{(j-N-1)!}{w(t)j!}dt.
\end{equation}

For $N\geq 1,$ it is then straightforward to see that
there is a constant, $M_{\Gamma,N},$ depending only on
$\Gamma$ and $N$ such that
\begin{equation}
  |e_N(z)|\leq M_{\Gamma,N}\|\varphi\|_{\cC^{N+1}(\Gamma)}|z|^{N+1}
\end{equation}
for any $z$ with $|z|\leq r.$
This proves
\begin{theorem}
\label{thm:1}
 Let $\Gamma$ be a smooth, bounded curve in $\bbC$ such
  that $B_r(0)\subset\Gamma^c,$ but $re^{i\theta_0}\in\Gamma\cap
  bB_r(0),$ where $\Gamma^c$ denotes the complement of $\Gamma\subset
  \mathbb R^2$.
  For each positive integer $N$ there is a constant
  $M_{\Gamma,N}$ such that, for $\varphi\in\cC^{N+1}(\Gamma),$ the error
in the truncated Taylor series approximation is given by
\begin{equation}\label{eqn13}
  \left|\frac{1}{2\pi i}\int\limits_{\Gamma}\frac{\varphi(w)dw}{w-re^{i\theta_0}}
-\sum_{j=0}^N\frac{ f^{(j)}_-(0)}{j!}r^je^{ij\theta_0}\right|\leq
  M_{\Gamma,N}\|\varphi\|_{\cC^{N+1}(\Gamma)}r^{N+1}.
\end{equation}
\end{theorem}

\noindent
Note that the coefficients of the expansion can be obtained by
evaluating the non-singular integrals:
\begin{equation}\label{eqn14}
  \frac{ f^{(j)}_-(0)}{j!}=\frac{1}{2\pi i}\int\limits_0^L
\frac{\varphi(w(t))w'(t)dt}{[w(t)]^{j+1}}.
\end{equation}
Error estimates for the numerical computation of these integrals are discussed in section
\ref{sec:quad}.

We turn now to a second approach for estimating the error in truncating the Taylor series,
which uses the Fourier expansion of $f_-$ on $bB_r(0)$ directly. As $f_-$ is
holomorphic in $B_r(0)$ its Fourier series on the boundary only has  terms
with non-negative exponents:
\[
  f_-(re^{i\theta})=\sum_{j=0}^{\infty}a_je^{ij\theta},
  \qquad
  \text{where}
  \qquad
  a_j=\frac{1}{2\pi}\int\limits_{0}^{2\pi}f_-(re^{i\theta})e^{-ij\theta}d\theta.
\]
We again want to show that
\begin{equation}\label{eqn17}
  \left|f_-(re^{i\theta})-\sum_{j=0}^{N}a_je^{ij\theta}\right|=O(r^{N+1}),
\end{equation}
which is evidently a matter of estimating the Fourier coefficients
themselves.

Integrating by parts, and assuming $j>0$, we have
\begin{equation}
  a_j=\frac{1}{2\pi }\left(\frac{r}{j}\right)
\int\limits_{0}^{2\pi}f^{(1)}_-(re^{i\theta})e^{-i(j-1)\theta}d\theta.
\end{equation}
Note that $e^{-ij\theta}$ has been replaced by
$e^{-i(j-1)\theta},$ which is a consequence of the fact that
\begin{equation}
  \pa_{\theta}f_-(re^{i\theta})=ire^{i\theta}f_-^{(1)}(re^{i\theta}).
\end{equation}
This shift in degree reappears when we consider spherical harmonic
expansions below.

If $f_-$ has $N+1$ derivatives on $bB_r(0),$ then
we can repeat this $(N+1)$-times to obtain that, for $j\geq N+1,$
\begin{equation}
  a_j=\frac{1}{2\pi i}\left(\frac{r^{N+1}(j-N-1)!}{j!}\right)
\int\limits_{0}^{2\pi}f^{(N+1)}_-(re^{i\theta})e^{-i(j-N-1)\theta}d\theta,
\end{equation}
from which the estimate in~\eqref{eqn17} is immediate. The implied
constant depends on $\|f_-\|_{\cC^{N+1}(B_r(0))},$ which we know from the
mapping properties of the Cauchy integral operator is of the order
$O(\|\varphi\|_{\cC^{N+1,\beta}(\Gamma)})$ for any $\beta>0,$ see Theorem 1.4
on page 147 of~\cite{GCWen}. We
therefore get an estimate slightly weaker than that obtained
in~\eqref{eqn13}.
\subsection{Harmonic Layer Potentials}

We now apply similar ideas to the evaluation of single and double
layer potentials on the boundaries of planar regions.  The fundamental
solution for the Laplace operator is given by
\begin{equation}
  G_0(z,w)=\frac{1}{2\pi}\log|z-w|=\frac{1}{2\pi}\Re\log(z-w).
\end{equation}
While the complex log is multi-valued,
its real part and its complex derivative are globally defined as single-valued
analytic functions.  As
before we let $D_-$ be a bounded region with boundary $\Gamma.$
Assuming that $\varphi$ is real valued, the single layer potential with
density supported on $\Gamma$ is defined as
\begin{equation}
  u(\xi_1,\xi_2)=\frac{1}{2\pi}\Re\left[\int\limits_0^L\varphi(w(t))\log(w(t)-z)dt\right],
\end{equation}
where $z = \xi_1+i\xi_2$ and $w(t)$ is again an arc-length parametrization of $\Gamma,$
viewed as lying in the complex plane.

We assume (without loss of generality) that the origin $0 \in D_-$,
which we take to be the center of our power series expansion.
If $B_r(0)\subset D_-,$ then for $z\in B_r(0)$ we can rewrite this
integral as
\begin{equation} \label{eqnusep}
  u(\xi_1,\xi_2)=\frac{1}{2\pi}\Re\left[\int\limits_0^L\varphi(w(t))
\log\left(1-\frac{z}{w(t)}\right)dt\right]+ A_0,
\end{equation}
where
\begin{equation}\label{eqn24}
  A_0=\frac{1}{2\pi}\int\limits_0^L\varphi(w(t))\log|w(t)|dt.
\end{equation}
Let us denote by $u_1$ the first term in~\eqref{eqnusep}.  As $|z/w(t)|<1,$ the
power series expansion for the principal branch of $\log(1-\zeta)$ about $\zeta=0,$ shows that $u_1$
has the following power series expansion:
\begin{equation}
  u_1(\xi_1,\xi_2)=\frac{-1}{2\pi}\Re\left[\int\limits_0^L\varphi(w(t))\sum_{j=1}^{\infty}
\frac{1}{j}\left(\frac{z}{w(t)}\right)^jdt\right].
\end{equation}
To approximate $u(\xi_1,\xi_2)$ for $z\in\overline{B_r(0)}$ we use $A_0$ plus a finite partial sum of this series, giving an error term:
\begin{multline} \label{harm2dtrunc}
  e_N(\xi_1,\xi_2)=u(\xi_1,\xi_2)-\left[A_0-\frac{1}{2\pi}\Re\left\{\int\limits_0^L\varphi(w(t))\sum_{j=1}^{N}
\frac{1}{j}\left(\frac{z}{w(t)}\right)^jdt\right\}\right]=\\
\frac{1}{2\pi}\sum_{j=N+1}^{\infty}\Re\left[\int\limits_0^L\varphi(w(t))
\frac{1}{j}\left(\frac{z}{w(t)}\right)^jdt\right]
\end{multline}
Using the identity in equation~\eqref{eqn11.01}
we can integrate by parts $N$-times to obtain:
\begin{equation}
  e_N(\xi_1,\xi_2)=\frac{1}{2\pi}\sum_{j=N+1}^{\infty}
\frac{z^{N+1}(j-N)!}{j!}\Re\left[\int\limits_0^LD_t^N[\varphi(w(t))]
\left(\frac{z}{w(t)}\right)^{j-N-1}\frac{dt}{w(t)}\right],
\end{equation}
where $D_t$ is again the  differential operator defined in~\eqref{eqn12.01}.
As $r\to 0,$
\begin{equation}
  \int\limits_{\Gamma}\frac{dt}{|w(t)|}\propto -\log r.
\end{equation}
The preceding estimates yield the following theorem.

\begin{theorem} \label{harm2dthm}  Suppose that $\Gamma$ is a smooth, bounded curve
  embedded in $\bbR^2,$ such that $B_r(0)\subset\Gamma^c,$ but
  $z_0 = re^{i\theta_0}\in bB_r(0)\cap\Gamma.$ For $N$ a positive integer,
  there is a constant $M_{N,\Gamma}$ so that if
  $\varphi\in\cC^{N}(\Gamma),$ then
\begin{multline}\label{eqn53.3}
  \left|\int\limits_{\Gamma}G_0(z_0,y)\varphi(y)ds(y)-
\left[A_0-\frac{1}{2\pi}\Re\left\{\int\limits_0^L\varphi(w(t))\sum_{j=1}^{N}
\frac{1}{j}\left(\frac{z}{w(t)}\right)^jdt\right\}\right]\right|\\
\leq M_{N,\Gamma}\|\varphi\|_{\cC^N(\Gamma)}r^{N+1}\log\frac{1}{r}.
\end{multline}
The constant $A_0$ is defined in~\eqref{eqn24}.
\end{theorem}

The double layer potential is defined as
 \begin{equation}\label{eqn32}
   v(\xi_1,\xi_2)=\frac{1}{2\pi}\int\limits_{\Gamma}\pa_{\bn_w}\log|w-z|\varphi(w(t))dt,
 \end{equation}
where $\bn$ is the outward normal along $bD_-.$
As $\log(w-z)$ is analytic, the Cauchy-Riemann equations imply that
\begin{equation}
  \pa_{\bn_w}\log|w-z|=\Re[\pa_{\bn_w}\log(w-z)]=-\Im[\pa_{\btau_w}\log(w-z)],
\end{equation}
where $\btau$ is the unit tangent vector along $bD_-.$ Thus, we may write
\begin{equation}
    \pa_{\bn_w}\log|w-z|=-\Im\frac{w'(t)}{w(t)-z}.
\end{equation}
It follows that the estimate proved in the previous section for
the Cauchy transform can be applied to analyze the error when replacing the
double layer potential, $v(\xi^0_1,\xi^0_2)$ with the finite partial sum:
\begin{equation}\label{eqn35}
v_{(N)}(r,\theta_0)= - \frac{1}{2\pi}\sum_{j=0}^N\Im\left[\int\limits_{\Gamma}\left(\frac{re^{i\theta_0}}{w(t)}
\right)^j\frac{\varphi(w(t))w'(t)dt}{w(t)}\right]
\end{equation}
where $(r,\theta_0)$ are the polar coordinates of $(\xi_1^0,\xi_2^0)$ with respect
to the origin (which is the center of the local expansion).

\begin{corollary}
 Let $\Gamma$ be a smooth, bounded curve in $\bbC$ such that $B_r(0)\subset\Gamma^c,$
with $z_0 = \xi^0_1+i\xi^0_2 = re^{i\theta_0}\in\Gamma\cap bB_r(0).$
For each positive integer $N$ there is a constant $M_{N,\Gamma}$ such that
if $\varphi\in\cC^{N+1}(\Gamma),$  then
\begin{equation}
  |v(\xi^0_1,\xi^0_2)-v_{(N)}(r,\theta_0)|\leq
M_{N,\Gamma}r^{N+1}\|\varphi\|_{\cC^{N+1}(\Gamma)},
\end{equation}
where $v$ is defined in~\eqref{eqn32} and $v_{(N)}$ is defined in~\eqref{eqn35}.
\end{corollary}

\subsection{Layer Potentials for the Helmholtz Equation}
\label{sec:two-d-helmholtz}

Our last example in two dimensions concerns layer potentials that satisfy the
Helmholtz equation:
\begin{equation}
  (\Delta +k^2)u=0.
\end{equation}
The fundamental solution, for $k\neq 0,$ is given by
the zeroth order Hankel function of the first kind:
\begin{equation}
  G_k(x,y)=\frac{i}{4} H^{(1)}_0(k|x-y|).
\end{equation}
It is well-known that the kernel $G_k(x,y)$ defines a classical
pseudodifferential operator of order $-2,$ see~\cite{Taylor2}.
For our purposes, the most important fact about this function is the Graf addition
theorem, see equation 9.1.79 in~\cite{AS}, which gives the series representation:
\begin{equation}\label{eqn40.01}
  H^{(1)}_0(k|x-y|)=\sum_{l=-\infty}^{\infty}J_{|l|}(k|x|)H^{(1)}_{|l|}(k|y|)e^{il(\theta-\phi)},
\end{equation}
valid so long as
$|x|<|y|\text{ where }x=|x|e^{i\theta}\text{ and }y=|y|e^{i\phi}.$

Here $J_l$ and $H_l^{(1)}$ denote the usual Bessel and Hankel functions of the
first kind, respectively.

We assume (without loss of generality) that the origin $0 \in D_-$, which we
take to be the center of our series expansion.  With these
assumptions about $D_-$ and $B_r(0),$ we define the single layer potential:
\begin{equation}\label{eqn41}
\begin{split}
  u(x)&=\int\limits_{\Gamma}\varphi(y)G_k(x,y)ds(y)\\
&=\sum_{l=-\infty}^{\infty}J_{|l|}(k|x|)e^{il\theta}\int\limits_{\Gamma}H^{(1)}_{|l|}(k|y|)
e^{-il\phi_y}\varphi(y)ds(y)\\
&=\sum_{l=-\infty}^{\infty}\alpha_lJ_{|l|}(k|x|)e^{il\theta}.
\end{split}
\end{equation}
The second equation follows from~\eqref{eqn40.01}
The series in the last line is valid for $x\in B_r(0),$ with the coefficients
$\{\alpha_l\}$ defined by the integrals:
\begin{equation}\label{eqn43}
  \alpha_l=\int\limits_{\Gamma}H^{(1)}_{|l|}(k|y|)
e^{-il\phi_y}\varphi(y)ds(y).
\end{equation}

We suppose that $\varphi\in\cC^{N,\beta}(\Gamma),$ for some $\beta>0,$
so that $u_{\pm}$ belongs to
$\cC^{N+1,\beta}(\overline{D_{\pm}}).$ If we assume that
$x_0=re^{i\theta_0}\in bB_r(0)\cap\Gamma,$ then, in the QBX approach,
we would like to use the partial sum of the series:
\begin{equation}
 \sum_{l=-N}^{N}\alpha_lJ_{|l|}(kr)e^{il\theta_0}
\end{equation}
as an approximation for $u_-(x_0).$
Because the Fourier series representation is unique, it follows that
\begin{equation}
  \alpha_lJ_{|l|}(kr)=\frac{1}{2\pi}\int\limits_{0}^{2\pi}
u_-(r\cos\theta,r\sin\theta)e^{-il\theta}.
\end{equation}
We now seek, as above, to estimate the error
\begin{equation} \label{helmtrunc}
  \left|u_-(x_0)-\sum_{l=-N}^{N}\alpha_lJ_{|l|}(kr)e^{il\theta_0}\right| .
\end{equation}

We make use of the integration by parts formula:
\begin{multline}
  \frac{1}{2\pi}\int\limits_{0}^{2\pi}
u_-(r\cos\theta,r\sin\theta)e^{-il\theta}=\\
\frac{r}{2\pi\cdot l}\int\limits_{0}^{2\pi}
\left([\pa_{\bz}u_-](r\cos\theta,r\sin\theta)e^{-i(l+1)\theta}-
[\pa_{z}u_-](r\cos\theta,r\sin\theta)e^{-i(l-1)\theta}\right)d\theta,
\end{multline}
where
\begin{equation}
  \pa_{z}=\frac{1}{2}(\pa_x+i\pa_y)\text{ and }
\pa_{\bz}=\frac{1}{2}(\pa_x-i\pa_y).
\end{equation}
If we integrate by parts, in this manner, $N+1$ times,
then the order $l$ term, where $N<|l|,$ produces an integrand of the form:
\begin{equation}
r^{N+1}\sum_{j=0}^{N+1} c_{l,N+1,j}
e^{i(l+2j-N-1)\theta}[\pa_z^j\pa_{\bz}^{N+1-j}u_-](r\cos\theta,r\sin\theta).
\end{equation}
Each coefficient $c_{l,N+1,j}$ is a sum of $\left(\begin{matrix}
    N+1\\j\end{matrix}\right)$ terms of the form:
\begin{equation}
  \frac{1}{l(l+\epsilon_1)\cdots(l+\epsilon_N)}\text{ where }
\epsilon_j\in\{-N,1-N,\dots,N-1,N\},
\end{equation}
with each  satisfying the upper bound:
\begin{equation}
  |c_{l,N+1,j}|\leq\frac{(|l|-N-1)!}{|l|!}.
\end{equation}
Thus, for $|l|>N,$ summing over $j$ gives the estimate
\begin{equation}
  |\alpha_lJ_{|l|}(kr)|\leq \frac{(2r)^{N+1}(|l|-N-1)!\cdot\|u_-\|_{\cC^{N+1}(\overline{B_r(0)})}}{|l|!}.
\end{equation}
Hence, for $N\geq 1,$ there is a constant $M_N$ so that
\begin{equation}\label{eqn52}
  \left|u_-(re^{i\theta_0})-\sum_{l=-N}^{N}\alpha_lJ_{|l|}(kr)e^{il\theta_0}\right|\leq
M_Nr^{N+1}\|u_-\|_{\cC^{N+1}(\overline{B_r(0)})}.
\end{equation}
The constant $M_N$ in this estimate does not depend on the frequency, but we
have not yet estimated the error in terms of the data $\varphi.$ The layer potential operator
corresponding to
$G_k$ maps $\cC^{N,\beta}(\Gamma)$ to $\cC^{N+1,\beta}(\overline{D_-})$ for any
$\beta>0,$ see~\cite{ColtonKress}. Hence the preceding estimates yield the
following result.
\begin{theorem}
\label{thm:3}
 Suppose that $\Gamma$ is a smooth, bounded curve
  embedded in $\bbR^2,$ such that $B_r(0)\subset\Gamma^c,$ with
  $re^{i\theta_0}\in bB_r(0)\cap\Gamma.$ For $k\in[0,\infty),$ $N$ a positive integer,
  and $\beta>0,$ there is a constant $M'_{N,\beta}(k)$ so that if
  $\varphi\in\cC^{N,\beta}(\Gamma),$ then
\begin{multline}\label{eqn53}
  \left|\int\limits_{\Gamma}G_k(re^{i\theta_0},y)\varphi(y)ds(y)-\sum_{l=-N}^{N}\alpha_lJ_{|l|}(kr)e^{il\theta_0}\right|\leq\\
M'_{N,\beta}(k)r^{N+1}\|\varphi\|_{\cC^{N,\beta}(\Gamma)}.
\end{multline}
Here the coefficients $\{\alpha_l\}$ are given by~\eqref{eqn43}.
\end{theorem}
\begin{remark}
  \label{rem:k-dependence-fourier}
  The dependence on the wave number $k$ of the constant in this
  estimate arises from the dependence on $k$ of the norm of the operator
  \begin{equation}\label{eqn2.52.001}
    G_k:\cC^{N,\beta}(\Gamma)\longrightarrow\cC^{N+1,\beta}(\overline{D_-}).
  \end{equation}
Let $m_{N,\beta}(k)$ denote this norm. One can show that there are constants
$\tm_{N,\beta}$ independent of $k$ so that
\begin{equation}
  m_{N,\beta}(k)\leq (1+k^{N+1})\tm_{N,\beta}.
\end{equation}
 As it is tangential to the main thrust
of this paper, we only explain briefly  how to obtain such an estimate  in the simplest
case, where $bD$ is the $x_1$-axis. First observe that $G_k$ is a convolution
operator
\begin{equation}
  G_k\varphi(x_1,x_2)=\int\limits_{-\infty}^{\infty}g_k(x_1-y_1,x_2)\varphi(y_1)dy_1.
\end{equation}
Using the equation, $(\pa_{x_1}^2+\pa_{x_2}^2)u=-k^2u,$ at the expense of
introducing powers of $k,$ we can replace normal derivatives on $u$ with
tangential derivations, which can then be shifted to the data. That is
\begin{equation}
\begin{split}
 (-1)^l \pa_{x_1}^m\pa_{x_2}^{2l}u(x_1,x_2)&=\pa_{x_1}^m(\pa_{x_1}^2+k^2)^lu(x_1,x_2)\\
&=\int\limits_{-\infty}^{\infty}g_k(x_1-y_1,x_2)((\pa_{y_1}^2+k^2)^l\pa_{y_1}^m\varphi(y_1)dy_1.
\end{split}
\end{equation}

As they are somewhat simpler, we give the estimates for 1-dimensional $L^2$-Sobolev spaces.
 The Fourier representation for the operator $G_k$ is
 \begin{multline}
   G_k\varphi(x_1,x_2)=\\
C_-\int\limits_{|\xi_1|<k}\frac{e^{ix_1\xi_1}e^{ix_2\sqrt{k^2-\xi_1^2}}\hvarphi(\xi_1)d\xi_1}{\sqrt{k^2-\xi_1^2}}
   +C_-\int\limits_{|\xi_1|>k}\frac{e^{ix_1\xi_1}e^{-x_2\sqrt{\xi_1^2-k^2}}\hvarphi(\xi_1)d\xi_1}{\sqrt{\xi_1^2-k^2}}.
 \end{multline}
Here the constants $C_-, C_+$ do not depend on $k.$ Using these two relations
we see easily that there are constants independent of $k$ so that, for $x_2\geq
0$ we have
\begin{equation}
  \int\limits_{-\infty}^{\infty}|\pa_{x_1}^m\pa_{x_2}^{2l+1}u(x_1,x_2)|^2dx_1
\leq C_{m,l}k^{2l}\|\varphi\|^2_{W^{m+2l,2}(\bbR)}.
\end{equation}
The dependence on $k$ of the norms of the operators in~\eqref{eqn2.52.001} is
quite similar.
\end{remark}

An estimate for the
double layer
\begin{equation}
  v(x)=\int\limits_{\Gamma}\pa_{\bn_y}G_k(x,y)\varphi(y)ds(y)
\end{equation}
is obtained just as easily.
In $D_-$, the double layer has an expansion of the form
\begin{equation}
  v_-(x)=\sum_{l=-\infty}^{\infty}\alpha_lJ_{|l|}(k|x|)e^{il\theta},
\end{equation}
valid for $x\in B_r(0),$ where now the coefficients
$\{\alpha_l\}$ are defined by the integrals:
\begin{equation}\label{eqn57}
  \alpha_l=\int\limits_{\Gamma}\pa_{\bn_y}[H^{(1)}_{|l|}(k|y|)
e^{-il\phi_y}]\varphi(y)ds(y).
\end{equation}
Arguing exactly as before we can show that, for $N\geq 2,$ we have the
estimate:
\begin{equation}
  \left| v(x_0)-
\sum_{l=1-N}^{N-1}\alpha_lJ_{|l|}(k|x|)e^{il\theta}\right|
\leq M_Nr^{N}\|v_-\|_{\cC^{N}(\overline{B_r(0)})}.
\end{equation}
Because the double layer defines an operator of order $-1$ we have
\begin{theorem}
\label{thm:4}
Suppose that $\Gamma$ is a smooth, bounded curve
  embedded in $\bbR^2,$ such that $B_r(0)\subset\Gamma^c,$ but
  $re^{i\theta_0}\in bB_r(0)\cap\Gamma.$ For $k\in[0,\infty),$ $N$ a
  positive integer, and $\beta>0,$ there is a constant
  $M''_{N,\beta}(k)$ so that if $\varphi\in\cC^{N,\beta}(\Gamma),$ then
\begin{multline}\label{eqn53d}
  \left|\int\limits_{\Gamma}\pa_{\bn_y}G_k(re^{i\theta_0},y)\varphi(y)ds(y)-\sum_{l=1-N}^{N-1}\alpha_lJ_{|l|}(kr)e^{il\theta_0}\right|\leq\\
M''_{N,\beta}(k)r^{N}\|\varphi\|_{\cC^{N,\beta}(\Gamma)}.
\end{multline}
Here the coefficients $\{\alpha_l\}$ are given by~\eqref{eqn57}.
\end{theorem}
\begin{remark}
  \label{rem:k-dependence-coefficient}
  The coefficient in this estimate also satisfies a bound of the form
  \begin{equation}
    M''_{N,\beta}(k)\leq m''_{N,\beta}(1+k^N),
  \end{equation}
where $m''_{N,\beta}$ is independent of $k.$
\end{remark}
\section{Three-dimensional Layer Potentials}
\label{sec:three-d}
In this section we consider the evaluation of layer potentials arising
from the Helmholtz equation in three dimensions. The ``outgoing'' fundamental
solution is given by
\begin{equation}
  G_k(x,y)=\frac{e^{ik|x-y|}}{4\pi |x-y|}
\end{equation}
for $\Im k\geq 0.$ The relevant formul{\ae} in the harmonic ($k=0$) case are somewhat
different, but follow the same lines as those for $k\neq 0.$ We leave that
analysis to the interested reader.

As in the two-dimensional case, this Green's function also
has a classical expansion in terms of spherical harmonics \cite{MF}:
\begin{equation}
  G_k(x,y)= ik \sum_{l=0}^\infty
j_{l}(k|x|)h_{l}(k|y|)\sum_{m=-l}^lY_l^m(\theta_1,\phi_1)
Y_l^{-m}(\theta_2,\phi_2)
\end{equation}
assuming $|x|<|y|$, with
\begin{equation}
  x=|x|\omega(\theta_1,\phi_1)\text{ and }
  y=|y|\omega(\theta_2,\phi_2) \, ,
\end{equation}
where
\begin{equation}
\omega(\theta,\phi)=(\sin\phi\cos\theta,\sin\phi\sin\theta,\cos\phi).
\end{equation}
Here, $j_{l}$ and $h_{l}$ denote the spherical Bessel and Hankel
functions of the first kind, respectively \cite{AS}.
The normalized spherical harmonics $Y_l^m$ are defined for $|m|\leq l$ by
\begin{equation}
  Y_l^m(\theta,\phi)=\sqrt{\frac{2l+1}{4\pi}\cdot\frac{(l-|m|)!}{(l+|m|)!}}
  P_l^{|m|}(\cos\phi)e^{im\theta},
\end{equation}
where $P_l^m$ is the associated Legendre function of degree $l$ and order $m$.

We let $D_-$ denote a bounded region in $\bbR^3$ with smooth boundary
$\Gamma.$ We assume, as in the two-dimensional case, that the origin $0 \in D_-$
and that $B_r(0)\subset D_-,$ with $x_0\in
bB_r(0)\cap\Gamma.$ If $\varphi$ is an integrable function defined on
$\Gamma,$ then a single layer potential with this density is given by
\begin{equation}\label{eqn67}
  \begin{split}
    f(x)&=\int\limits_{\Gamma}G_k(x,y)\varphi(y)dS(y)\\
&=\sum_{l=0}^{\infty} j_{l}(k|x|)
\sum_{m=-l}^l\alpha_{lm}Y_l^m(\theta_1,\phi_1).
  \end{split}
\end{equation}
The series representation is valid in $B_r(0),$ with
\begin{equation}\label{eqn70}
  \alpha_{lm}=ik \int\limits_{\Gamma}
h_{l}(k|y|)\varphi(y)
Y_l^{-m}(\theta_y,\phi_y)dS(y) .
\end{equation}
Once again, we would like to approximate $f(x_0)$
by the partial sum of this expansion:
\begin{equation}\label{eqn71.0}
   f(x_0)\approx f^{(N)}(x_0)=\sum_{l=0}^{N} j_{l}(k|x_0|)
\sum_{m=-l}^l\alpha_{lm}Y_l^m(\theta_0,\phi_0).
\end{equation}
To estimate the error $|f(x_0)- f^{(N)}(x_0)|$, we use an argument
very much like that used to obtain~\eqref{eqn52} and~\eqref{eqn53}.
Before we turn to that task, however, let us compute the series expansion
of the double layer
\begin{equation}
  v(x)=\int\limits_{\Gamma}\pa_{\bn_y}G_k(x,y)\varphi(y)dS(y).
\end{equation}
$v(x)$  also has an expansion like that in the second line of~\eqref{eqn67},
\begin{equation}\label{eqn71}
  v(x)=\sum_{l=0}^{\infty} j_{l}(k|x|)
\sum_{m=-l}^l\alpha_{lm}Y_l^m(\theta_1,\phi_1),\text{ for }|x|<r,
\end{equation}
where the coefficients are now defined by the integrals:
\begin{equation}
  \alpha_{lm}=ik \int\limits_{\Gamma}\pa_{\bn_y}\left[
h_{l}(k|y|) Y_l^{-m}(\theta_y,\phi_y)\right]\varphi(y)dS(y).
\end{equation}
The value of $v(x_0)$ can again be approximated by taking a partial
sum of the series in~\eqref{eqn71}.

The error estimates for both of these series approximations come from
determining the $r$ dependence of the coefficients in these expansion.
We carry out the analysis for the single layer and leave the details for the
double layer to the interested reader.

The necessary estimates follow from integration by parts using the fact that:
\begin{equation}\label{eqn74}
  L_{\pm}Y_l^m=\sqrt{(l\mp m)(l\pm m+1)}Y_l^{m\pm 1},
\end{equation}
where the vector fields are defined in $(\theta,\phi)$--coordinates by:
\begin{equation}
  L_{\pm}=e^{\pm i\theta}\left[\pa_{\phi}\pm i\cot\phi\pa_{\theta}\right].
\end{equation}
An elementary calculation shows that, as operators on $L^2(S^2),$ the
adjoints of $L_{\pm}$ are $-L_{\pm}$. In integrating by parts, we
need to apply these vector fields to $g(r\omega(\theta,\phi)),$ where
$g$ is a function smooth in the closure of the sphere of radius $r.$
An elementary calculation gives:
\begin{equation}
  L_{\pm}g(r\omega(\theta,\phi))=r\cos\phi[g_x\pm ig_y](r\omega(\theta,\phi))-
re^{\pm i\theta}\sin\phi g_z(r\omega(\theta,\phi)).
\end{equation}
It is useful to observe that:
\begin{equation}
  \cos\phi= \sqrt{\frac{4\pi}{3}}Y_1^0\text{ and }
e^{\pm i\theta}\sin\phi=\sqrt{\frac{8\pi}{3}}Y_1^{\pm 1},
\end{equation}
and therefore
\begin{equation}\label{eqn78}
  L_{\pm}g(r\omega(\theta,\phi))=r\sqrt{\frac{4\pi}{3}}[Y_1^0(g_x\pm ig_y)-
\sqrt{2}Y_1^{\pm 1} g_z](r\omega(\theta,\phi)).
\end{equation}

The Clebsch-Gordan relations \cite{ROSE} (or \cite{dlmf_2010},
\dlmf{34.3.20}, \dlmf{34.2.1}, \dlmf{34.2.3}), show that there are
collections of coefficients $\{c^+_{jlm;n}; c^-_{jlm;n}\}$ so that for
$j\in\{-1,0,1\}$ we have the identities
\begin{equation}\label{eqn79}
  Y_1^jY_l^{-m}=\sum_{n=1-l}^{l-1}c^-_{jlm;n}Y_{l-1}^{n}+\sum_{n=-(1+l)}^{1+l}c^+_{jlm;n}Y_{l+1}^{n}.
\end{equation}
An elementary computation shows that
\begin{equation}
  |Y_1^jY_l^{-m}(\omega)|^2\leq\frac{3}{4\pi(1+|j|)}|Y_l^{-m}(\omega)|^2.
\end{equation}
Integrating this inequality over the sphere, and using the
orthogonality relations satisfied by the functions $\{Y_l^m\},$ we
easily deduce that for $j\in\{-1,0,1\},$
\begin{equation}\label{eqn81}
 \sum_{n=1-l}^{l-1}|c^-_{jlm;n}|^2+
\sum_{n=-(l+1)}^{l+1}|c^+_{jlm;n}|^2\leq  \frac{3}{4\pi(1+|j|)}.
\end{equation}

As before, we observe that the representation
$f(r\omega(\theta,\phi))$
of the single layer potential in~\eqref{eqn67}
is unique, and therefore the coefficients
can be computed by integration on a sphere of radius $r$ about the origin:
\begin{equation}\label{eqn83}
  \talpha_{lm}(r)=
j_{l}(kr) \alpha_{lm}=\int\limits_{0}^{2\pi}\int\limits_{0}^{\pi}
  f(r\omega(\theta,\phi))Y^{-m}_{l}(\theta,\phi)\sin\phi d\phi d\theta.
\end{equation}
For simplicity assume that $m\geq 0,$ and use~\eqref{eqn74} to see that
\begin{equation}
  L_+Y_l^{-m}=\sqrt{(l+m)(l-m+1)}Y_l^{1-m}.
\end{equation}
Integrating by parts in the integral defining $\talpha_{lm}(r)$ we see that
\begin{equation}
  \talpha_{lm}(r)=
  -\int\limits_{0}^{2\pi}\int\limits_{0}^{\pi}
  L_+[f(r\omega(\theta,\phi))]\frac{Y_l^{1-m}(\theta,\phi)}{\sqrt{(l+m)(l-m+1)}}
dS(\theta,\phi)
\end{equation}
Using the relations in~\eqref{eqn78} and~\eqref{eqn79},
we see that the integrand becomes
\begin{multline}\label{eqn85}
  r \sqrt{\frac{4\pi}{3}}\frac{[(f_x+if_y)Y_1^0-\sqrt{2}f_z Y_1^1]Y_l^{1-m}}
{\sqrt{(l+m)(l-m+1)}}=\\
r \sqrt{\frac{4\pi}{3(l+m)(l-m+1)}}\Bigg\{
(f_x+if_y)
\left[\sum_{n=1-l}^{l-1}c^-_{0lm;n}Y_{l-1}^{n}+\sum_{n=-(1+l)}^{1+l}c^+_{0lm;n}Y_{l+1}^{n}\right]-\\
\sqrt{2}f_z
\left[\sum_{n=1-l}^{l-1}c^-_{1lm;n}Y_{l-1}^{n}+\sum_{n=-(1+l)}^{1+l}c^+_{1lm;n}Y_{l+1}^{n}\right]
\Bigg\}.
\end{multline}

The formula in~\eqref{eqn85} has several  notable features:
\begin{enumerate}
\item The entire quantity is multiplied by $r.$
\item The terms in these sums are of exactly the same types as those that
  appear in the original definition, with the integrals computing
  spherical harmonic expansion coefficients of derivatives of $f$
  along $|x|=r.$
\item Starting with a spherical harmonic of degree $l$ we obtain terms of
  degree $l-1$ and $l+1,$ showing that this procedure can be
  repeated $l$ times.
\item The coefficients $\{c^{\pm}_{jlm;n}\}$ can be bounded by the
  identity in~\eqref{eqn81}.
\end{enumerate}
Using these observations, we carry out integration by parts
$N+1$ times to obtain an estimate for the remainder:
\begin{equation}
  e^{(N)}(x_0)= f(x_0)- f^{(N)}(x_0)=\sum_{l=N+1}^{\infty}
\sum_{m=-l}^l\talpha_{lm}(r)Y_l^m(\theta_0,\phi_0) .
\end{equation}

Integrating the right hand side of~\eqref{eqn85} over the unit sphere leads to
the computation of certain spherical harmonic coefficients of $(f_x+if_y),$ and
$\sqrt{2}f_z$ restricted to the sphere of radius $r.$ In order to estimate the
coefficients $\{\talpha_{lm}(r)\},$ we need to introduce some notation.  We let
\begin{equation}
  D_1=\pa_x+i\pa_y,\,
 D_2=\pa_x-i\pa_y,\,
D_3=\sqrt{2}\pa_z.
\end{equation}
For a 3-multi-index $\gamma=(\gamma_1,\gamma_2,\gamma_3),$ we let
$\{\halpha_{lm}^{(\gamma)}(r)\}$ denote the spherical harmonic
coefficients of $D^\gamma f\restrictedto_{bB_r(0)},$ where $D^{\gamma}$
is the constant coefficient differential operator
$D_1^{\gamma_1}D_2^{\gamma_2}D_3^{\gamma_3},$ so that
\begin{equation}
  \halpha_{lm}^{(\gamma)}(r)=\int\limits_{0}^{2\pi}\int\limits_{0}^{\pi}
  [D^{\gamma}f](r\omega(\theta,\phi)) Y_l^m(\theta,\phi) \sin\phi d\phi d\theta.
\end{equation}

To obtain estimates for $\{\talpha_{lm}(r)\}$ that are proportional to
$r^{N+1},$ we need to integrate by parts in~\eqref{eqn83} $N+1$
times. Each integration by parts increases the number of terms by a
factor of $4.$ Whether one chooses $L_+$ or $L_-$ at each step is not
too important, but for the fact that for any pair $(l,m)$ one can
always choose the sign so that coefficient on the right hand side
of~\eqref{eqn74} does not vanish.

Our formula expresses
\begin{equation}
 \sum_{m=-l}^l \talpha_{lm}(r) Y_l^m(\theta,\phi)
\end{equation}
as a sum of $(2l+1)\times 4^{N+1}$ sums, each consisting of terms of the form
\begin{equation}\label{eqn90}
  \left(\frac{4\pi}{3}\right)^{\frac{N+1}{2}}
r^{N+1}\cdot\halpha^{(\gamma)}_{l_1v_1}(r)C_1\cdots C_{N+1}Y_l^{u_{N+1}}(\theta,\phi).
\end{equation}
Each matrix $C_q$ is
of the form
\begin{equation}
  C_q=\frac{c^{\pm}_{j_ql_qu_q;v_q}}{d_q}.
\end{equation}
For each $q,$ the pair $u_q,v_q$ are the matrix indices in $C_q,$
\begin{align}
j_q\in &\{-1,0,1\},\\
l-(N+1)\leq &l_q\leq l+(N+1), \text{ and therefore }\\
u_q\in\{-l_q,\dots, l_q\}&\text{ and }v_q\in\{-(l_q\pm 1),\dots, (l_q\pm 1)\}.
\end{align}
We think of $\halpha_{l_1v_1}^{(\gamma)}(r)$ as a $1\times [2(l_1\pm 1)+1]$
vector (indexed by $v_1$), and $Y_l^{u_{N+1}}$ as a $(2l+1)\times 1$ vector of
functions (indexed by $u_{N+1}$).  Moreover $\gamma$ is a
3-multi-index satisfying $|\gamma|=N+1$.

Finally, the denominator is of the form
\begin{equation}
  d_q=\sqrt{(l_q\mp u_q)(l_q\pm u_q+1)},
\end{equation}
with $+$ or $-$ chosen so that $d_q$ never vanishes. With this
understood, $d_q$ is easily seen to satisfy the estimate:
\begin{equation}\label{eqn93}
  \sqrt{2l_q}\leq d_q.
\end{equation}

Using~\eqref{eqn81} and ~\eqref{eqn93}, we can bound the Frobenius
norm of each matrix $C_q,$ in~\eqref{eqn90} (with $u_q$ and $v_q$ as the
matrix indices) by $\sqrt{3/4\pi}.$ Hence the absolute value of the quantity
in~\eqref{eqn90} is bounded by
\begin{equation}
  r^{N+1}\|\halpha^{(\gamma)}_{l_1\cdot}(r)\|\sqrt{\frac{2l+1}{4\pi}},
\end{equation}
where we use the fact that pointwise
\begin{equation}
  \sum_{m=-l}^l|Y_l^m(\theta,\phi)|^2=\frac{2l+1}{4\pi}.
\end{equation}

It therefore follows that for some pairs $\{(\gamma_j,p_j):\: j=1,\dots,
(2l+1)\cdot 4^{N+1}\}$ where $|\gamma_j|=N+1$ and $p_j\in\{-(N+1),\dots,
(N+1)\},$ chosen as described above, we have
\begin{equation}
  \left|\sum_{m=-l}^l \talpha_{lm}(r) Y_l^m(\theta,\phi)\right|\leq
\sqrt{\frac{2l+1}{4\pi}}r^{N+1}
\sum_{j=1}^{(2l+1)\cdot 4^{N+1}}\|\halpha^{(\gamma_j)}_{(l-p_j)\cdot}\|.
\end{equation}
Finally, collecting terms we see that there is a constant $K_N$ so
that:
\begin{equation}
  |e^{(N)}(x_0)|\leq K_Nr^{N+1}\sum_{\gamma\in J_{3,N+1}}\sum_{l=0}^{\infty}
\left[\sum_{m=-l}^l|\halpha^{(\gamma)}_{lm}(r)|^2\right]^{\frac 12}(l+N+1)^{\frac 32}
\end{equation}
Here $J_{3,N+1}$ is the collection of 3-multi-indices with
$|\gamma|=N+1$.
In order for these infinite sums to be finite, we need to assume that,
for each $\gamma\in J_{3,N+1}$ the sum satisfies:
\begin{equation}\label{eqn98}
  \sum_{m=-l}^l|\halpha^{(\gamma)}_{lm}(r)|^2= O\left(\frac{1}{l^{5+\delta}}\right),
\end{equation}
for some $\delta >0.$

The estimate in~\eqref{eqn98} holds if, for each $\gamma\in
J_{3,N+1},$ the function $D^{\gamma}f\restrictedto_{bB_r(0)}$
belongs to the $L^2$-Sobolev space $W^{3+\delta,2}(bB_r(0)).$ For a
single layer, this would follow from assuming that $\varphi\in
W^{3+N+\delta,2}(\Gamma).$ This represents a small loss of regularity
over the estimates we obtained in two dimensions. This loss would
appear to be a result of the added complexity of multiplying spherical
harmonics rather than exponentials. A more careful estimate of the
terms in~\eqref{eqn90} might provide a somewhat better result.
\begin{theorem}
\label{thm:5}
Let $k$ be a complex number with $\Im k\geq 0$ and
  let $\Gamma$ denote a smooth hypersurface in $\bbR^3$, such that
  $B_r(0)\subset\Gamma^c$ with $x_0=r\omega(\theta_0,\phi_0)
\in\Gamma\cap bB_r(0).$ For each
  positive integer $N$ and $\delta>0,$ there is a constant
  $M_{N,\delta}(k)$ such that
\begin{multline}\label{eqn99}
  \left|\int\limits_{\Gamma}G_k(x_0,y)\varphi(y)dS(y)-
    \sum_{l=0}^{N} j_{l}(kr)
    \sum_{m=-l}^l\alpha_{lm}Y_l^m(\theta_0,\phi_0)\right|\leq\\
  M_{N,\delta}(k)r^{N+1}\|\varphi\|_{W^{3+N+\delta,2}(\Gamma)},
\end{multline}
 so long as $\varphi\in W^{3+N+\delta,2}(\Gamma).$
The coefficients $\{\alpha_{lm}\}$ are given by~\eqref{eqn70}.
\end{theorem}
\begin{remark}
  \label{rem:k-dependence-3d}
  As before, the implied constants in the estimates on the
  coefficients $\{\halpha^{(\gamma)}_{lm}(r)\}$ in terms of $f$ do not depend
  on $k.$ The reason for this is that all the integrations-by-parts used to
  derive these estimates involve vector fields tangent to the sphere $S_r(0).$
  The dependence on $k$ of the constant $M_{N,\delta}(k)$ derives from the
  dependence on $k$ of the norm of the operator
  \begin{equation}
    G_k:W^{N+3+\delta,2}(\Gamma)\longrightarrow W^{N+4+\delta,2}(bB_r(0)).
  \end{equation}
Arguing as above we can show that there are constants $\tm_{N,\delta}$
independent of $k$ so that
\begin{equation}
  M_{N,\delta}(k)\leq \tm_{N,\delta}(1+k^{N+3}).
\end{equation}
\end{remark}

%
\section{Quadrature error}
\label{sec:quad}
The prior estimates establish bounds on the tails of the expansions used by
QBX, and hence on the error incurred in their truncation.  They do not take
into account that the coefficients in these expansions must still be computed
numerically. In this section, we derive an error estimate for this quadrature
problem. This estimate then helps to establish an upper bound on the ``mesh
spacing'' $h$ needed, leading to straightforward control of the overall error
in QBX, which is bounded above by the sum of these truncation and quadrature
errors.

Concentrating for the moment on the single layer for the
Laplace equation in two dimensions, let us first rewrite \eqref{harm2dtrunc} as
\begin{equation} \label{harm2dtrunc2}
  e_N(\xi_1,\xi_2)=u(\xi_1,\xi_2)-A_0-\Re\left\{\sum_{j=1}^{N}
a_j z^j \right\}= \Re\left\{
\sum_{j=N+1}^{\infty} a_j z^j \right\}
\end{equation}
where
\begin{align*}
  A_0&=\frac{1}{2\pi}\int\limits_0^L\varphi(w(t))\log|w(t)|dt,\\
  a_j &= \frac{1}{2\pi} \int\limits_0^L \frac{\varphi(w(t))}{j[w(t)]^j}dt,
\end{align*}
and $\Gamma=w([0,L))$.
Recall, also, that the center $z_c$ here is assumed to be the origin, for simplicity
of notation. In practice, what we wish to estimate is actually
\begin{equation} \label{harm2dtrunc3}
  e_N(\xi_1,\xi_2)= \left| u(\xi_1,\xi_2)-A^Q_0-\Re\left\{\sum_{j=1}^{N}
a^Q_j z^j \right\} \right|,
\end{equation}
where ${A^Q}_0$ and $a_j^Q$ are the {\em computed approximants} of $A_0$ and $a_j$.
We write this error as
\begin{multline}
e_N(\xi_1,\xi_2) =
  \left|
\left( u(\xi_1,\xi_2) - A_0-\Re\left\{\sum_{j=1}^{N} a_j z^j \right\}  \right) \right. \\
- \left. (A^Q_0- A_0)-\Re\left\{\sum_{j=1}^{N}[a^Q_j-a_j] z^j \right\}
\right| .
\label{error_anal}
\end{multline}

The first term is $O(r^{N+1})$ as shown in Theorem \ref{harm2dthm}.
The second term depends strongly on the specific quadrature rule used.
To facilitate straightforward generalization to complicated
geometries in more than two dimensions, we choose a composite
Gauss quadrature for this role.
If the boundary $\Gamma$ were divided into $M$ equal subintervals
of size $h$, and a $q$th order Gauss quadrature rule were used on each subinterval,
then
\begin{equation}
\label{quad:error}
\left| (A^Q_0- A_0)-\Re\left\{\sum_{j=1}^{N}[a^Q_j-a_j] z^j \right\} \right|
= C_q(N,\Gamma) \left( \frac{h}{4r} \right)^{2q} \|\varphi\|_{\cC^{2q}}.
\end{equation}
To see this, we note that on a curve segment $\Gamma_i$ of length $h$, the
standard estimate for $q$-point Gauss-Legendre quadrature
\cite{davis} yields
\begin{eqnarray}
\left| \int_{\Gamma_i} \frac{\varphi(w(t))}{j[w(t)]^j}dt -
 \sum_{n=1}^q
   \frac{\varphi(w(t_n))}{j[w(t_n)]^j} \, w_n \right|
\hspace{2in} \hfill
\nonumber \\
\hspace{2in} \leq
 \frac{h^{2q+1}}{j(2q+1)}\frac{(q!)^4}{(2q)!^3}
 \left| D^{2q} \frac{\varphi(w(t))}{[w(t)]^j} \right|_\infty, \hspace{.2in}
\label{GLestimate}
\end{eqnarray}
where $D^{q}$ denotes the $q^{\rm th}$ derivative of the integrand with respect to the
integration parameter along $\Gamma_i$, and $\{t_n,w_n\}$ denote the Gauss
nodes and weights scaled to $\Gamma_i$. From
Stirling's approximation
\[ \sqrt{2\pi} n^{n+\frac{1}{2}} e^{-n} < n! <
 2 \sqrt{\pi} n^{n+\frac{1}{2}} e^{-n}, \]
we have
\[
\left| \int_{\Gamma_i} \frac{\varphi(w(t))}{j[w(t)]^j}dt -
 \sum_{n=1}^q
   \frac{\varphi(w(t_n))}{j[w(t_n)]^j} \, w_n \right|
\leq \frac{h^{2q+1}}{(2q)! 4^{2q}}
 \left| D^{2q} \frac{\varphi(w(t))}{[w(t)]^j} \right|_\infty .
\]
Omitting a detailed computation of the derivative, we note that it satisfies
a bound of the general form
\[
 \left| D^{2q} \frac{\varphi(w(t))}{[w(t)]^j} \right|_\infty \leq
C_q(N,\Gamma) \frac{(2q)!}{r^{j+2q}}\|\varphi\|_{\cC^{2q}}.
\]
The constant in this estimate depends on the order of the local expansion $N$
and on the smoothness of the curve $\Gamma,$ through a polynomial in
$(w,w^{(1)},\dots,w^{(2q)}).$ Since the dependence of the constant on $N,$ and
the parametrization are rather complicated, but independent of $\varphi,$ we
leave it in the form above.

Adding up the errors from all subintervals yields
\eqref{quad:error}.
The total error incurred by the QBX method can therefore be written in the form
\[ E = O(r^{N+1}) + C_q(N,\Gamma)\left( \frac{h}{4r} \right)^{2q}\|\varphi\|_{\cC^{2q}}. \]
Fixing the order of the expansion $N$ and letting $r = 4h$ yields an error of the
form:
\begin{equation}
  E = O(h^{N+1}) + C_q(N,\Gamma) \left( \frac{1}{16} \right)^{2q}\|\varphi\|_{\cC^{2q}}.
  \label{laperrest}
\end{equation}
Implemented in this manner, the QBX method is asymptotic,
converging like a method of order $N+1$, until the error is dominated by the second term.
To obtain a {\em convergent} scheme one can, for example, set $r = \sqrt{h}$ so that
\[ E = O(h^{\frac{N+1}{2}}) +  C_q(N,\Gamma)\left( \frac{h}{16} \right)^{q}\|\varphi\|_{\cC^{2q}}. \]

\begin{remark}
Similar estimates can be obtained for
Helmholtz layer potentials and for three-dimensional cases.
From Remarks \ref{rem:k-dependence-fourier}, \ref{rem:k-dependence-coefficient} and \ref{rem:k-dependence-3d}, however, for large $k$, we would have to let
$kr = 4h$ to yield an error of the type (\ref{laperrest}).
This is not surprising since it is the product $kr$ that is a dimensionless
quantity in wave scattering.
\end{remark}

In practice, when dealing with complicated boundaries,
adaptive discretization is generally required, and the assumption that the boundary
is divided into equal subintervals must be relaxed. This does not change
the error analysis in a substantial way, and we refer the
reader to \cite{QBX} for details and examples.
\section{Conclusions}

It is well-known that local expansions of layer potentials are analytic away from the
boundary. There are surprisingly few results in the literature, however,
concerning the behavior of such expansions when evaluated
at the limit of their radii of convergence - namely at the nearest boundary point itself.
The estimates in Theorems
\ref{thm:1},\ref{harm2dthm}, \ref{thm:3}, \ref{thm:4}, and \ref{thm:5}
show that, despite the absence of a geometric decay parameter,
high order accuracy can be obtained in computing
layer potentials in a manner that is controlled by the smoothness of the boundary data.
The analysis presented here serves as a complement to
\cite{QBX}, providing a convergence theory for the quadrature method QBX.
We note that sharper results can be obtained \cite{barnett_2012}, if
some assumptions are made about analyticity of the data and the location of
the nearest singularity.
Finally, estimates of the type obtained here may be of interest in areas of
mathematics that involve smooth continuation, microlocal analysis, and potential
theory.

\section*{Acknowledgment}

The authors would like to thank Alex Barnett, Zydrunas Gimbutas and
Michael O'Neil for many useful discussions.

{\bibliographystyle{siam} {\bibliography{qbx2}}}

\end{document}